\pdfoutput=1
\RequirePackage{ifpdf}
\ifpdf 
\documentclass[pdftex]{sigma}
\else
\documentclass{sigma}
\fi

\numberwithin{equation}{section}

\newtheorem{Theorem}{Theorem}[section]
 \newtheorem{Corollary}[Theorem]{Corollary}
 \newtheorem{Lemma}[Theorem]{Lemma}
 \newtheorem{Proposition}[Theorem]{Proposition}

{ \theoremstyle{definition}
 \newtheorem{Definition}[Theorem]{Definition}
 \newtheorem{Example}[Theorem]{Example}
 \newtheorem{Examples}[Theorem]{Examples}
 
}

\usepackage{amsxtra}
\usepackage{tikz}
\usepackage{tikz-cd}

\DeclareMathOperator{\rank}{rank}

\DeclareMathOperator{\Adjoint}{Ad}

\DeclareMathOperator{\dimension}{dim}

\DeclareMathOperator{\automorphism}{Aut}

\newcommand\integral{\int}

\renewcommand{\le}{\leqslant}
\renewcommand{\ge}{\geqslant}

\begin{document}
\allowdisplaybreaks

\newcommand{\arXivNumber}{1702.02717}

\renewcommand{\PaperNumber}{029}

\FirstPageHeading

\ShortArticleName{A Characterisation of Smooth Maps into a Homogeneous Space}

\ArticleName{A Characterisation of Smooth Maps\\ into a Homogeneous Space}

\Author{Anthony D.~BLAOM}

\AuthorNameForHeading{A.D.~Blaom}

\Address{University of Auckland, New Zealand}
\Email{\href{mailto:anthony.blaom@gmail.com}{anthony.blaom@gmail.com}}
\URLaddress{\url{https://ablaom.github.io}}

\ArticleDates{Received June 25, 2021, in final form April 04, 2022; Published online April 10, 2022}

\Abstract{We generalize Cartan's logarithmic derivative of a smooth map from a manifold into a Lie group $G$ to smooth maps into a~homogeneous space $M=G/H$, and determine the global monodromy obstruction to reconstructing such maps from infinitesimal data. The logarithmic derivative of the embedding of a submanifold $\Sigma \subset M$ becomes an invariant of~$\Sigma $ under symmetries of the ``Klein geometry'' $M$ whose analysis is taken up in [{\it SIGMA} {\bf 14} (2018), 062, 36 pages, arXiv:1703.03851].}


\Classification{53C99; 22A99; 53D17}

\section{Introduction}
According to a theorem of \'Elie Cartan, a smooth map $f \colon \Sigma \rightarrow G$, from a
connected manifold $\Sigma$ into a Lie group $G$, is uniquely
determined by its {\it logarithmic derivative}, up to right
translations in $G$. This derivative, also known as the Darboux
derivative of $f$, is a one-form $\delta f$ on~$\Sigma $ taking values
in the Lie algebra of ${\mathfrak g} $ of $G$. Here we formulate and
prove a generalization of this result, Theorem \ref{summaryT}, to
smooth maps $f \colon \Sigma \rightarrow G/H$ into an arbitrary
homogeneous space $G/H$. Our generalization describes explicitly the
global obstruction to reconstructing such maps from infinitesimal
data, data that generalizes logarithmic derivatives (generalized
Maurer--Cartan forms).

In this introduction we generalize the notion of Maurer--Cartan forms
and their monodromy, and state the main existence Theorem~\ref{outlineT}. The proof is straightforward, apart from a question
about Lie algebroid integrability, which is addressed in Section
\ref{newsec}, by applying \cite{Crainic_Fernandes_03}. Uniqueness, up
to symmetry, is guaranteed under a mild topological condition on
$G/H$, but we must take some care to qualify what is meant by
``symmetry'', a task postponed to Section~\ref{infch}.

Cartan's theorem is commonly associated with his method of moving
frames for studying subgeometry. While the moving frames method can
be reinterpreted within the present framework, it is possible to study
subgeometry using the new theory without fixing frames or local
coordinates. Both frame and frame-free illustrations are given in a
sequel article \cite{Blaom_G}. It is instructive to review Cartan's
approach here. For more detail we recommend \cite{Ivey_Landsberg_03}.

\subsection*{Cartan's method of moving frames}

To classify, with a unified approach, the submanifolds of Euclidean
space, affine space, conformal spheres, projective space, and so on,
the ambient space $M$ is viewed as a homogeneous space~$G/H$, i.e., as
a ``Klein geometry''. Here $G$ is the group of symmetries of the
geometric structure on $M$, which acts transitively by assumption.

Using the group structure, one tries to replace the embedding $f
\colon \Sigma \rightarrow G/H$ of a sub\-ma\-ni\-fold~$\Sigma $ with certain
data defined just on~$\Sigma $ and amounting to an
infinitesimalization of the map~$f$. The infinitesimal data consists
of {\it invariants} of $f$~--- that is, the data depends on $f$ only
up to symmetries of $G/H$ ($G$-translations). However, these
invariants ought to be {\it complete}, in~the sense that they are
sufficient for the reconstruction of $f$, up to symmetry.

Cartan's method for finding a complete set of invariants is in two
steps. In~the first step one attempts to lift the embedding $f \colon
\Sigma \rightarrow M \cong G/H$ to a smooth map $\tilde f \colon
\Sigma \rightarrow G$:
\[
 \begin{tikzcd} {}& G \arrow{d}\\ \Sigma \arrow{ur}{\tilde f}
\arrow{r}{f} & G/H.
 \end{tikzcd}
\]
The lift, which is not unique, should be as canonical as possible,
to make the identification of invariants easier later on. For
example, given a curve in Euclidean three-space $f \colon [0,1]
\rightarrow {\mathbb R}^3$, one obtains a lift $\tilde f \colon [0,1]
\rightarrow G$ into the group of rigid motions by declaring $\tilde
f(t)$ to be the rigid motion mapping the Frenet frame of the curve at
$f(0)$ to the Frenet frame at $f(t)$~--- the ``moving frame''.

Now the basic infinitesimal invariant of a Lie group $G$ is the
Maurer--Cartan form, a one-form on $G$ taking values in its Lie algebra
${\mathfrak g} $. In~the second step of Cartan's procedure, one pulls
the Maurer--Cartan form back from $G$ to a one-form $\delta\tilde f $
on~$\Sigma $ using the lifted map $\tilde f \colon \Sigma \rightarrow
G$. By~Cartan's theorem recalled below, one can reconstruct $\tilde f
\colon \Sigma \rightarrow G$, and hence the map $f \colon \Sigma
\rightarrow G/H$, from a knowledge of $\delta\tilde f $ alone, which
accordingly encodes (indirectly) complete invariants for the
embedding.

\subsection*{Smooth maps into a Lie group}

Fix a Lie group $G$ and let ${\mathfrak g}$ denote its Lie algebra. A
{\it Maurer--Cartan form} on a smooth mani\-fold~$\Sigma $ is a
${\mathfrak g}$-valued one-form $\omega $ satisfying the Maurer--Cartan
equations,
\begin{equation*}
{\rm d} \omega_k + \sum_{i<j}c_k^{ij}\omega_i \wedge
\omega_j =0,
\end{equation*}
where the $\omega_k$ are components of $\omega$ with
respect to some basis of ${\mathfrak g} $, and $c_k^{ij}$ the
corresponding structure constants. We have written the Maurer--Cartan
equations as they are most commonly recognized, although this is not
best representation from the present point of view, as we shall~see.

The Lie group $G$ itself supports a unique right-invariant
Maurer--Cartan form $\omega_G$ that is the identity on
$T_e{G}={\mathfrak g} $. Every smooth map $f \colon \Sigma \rightarrow
G$ pulls $\omega_G$ back to a Maurer--Cartan form on~$\Sigma $, here
denoted $\delta f$. Since
\begin{equation*}
\delta f \bigg(\frac{\rm d}{{\rm d}t}x(t)\bigg)
=\frac{\rm d}{{\rm d}\tau}f(x(\tau))f(x(t))^{-1}\bigg|_{\tau=t},
\end{equation*}
or $\delta f = {\rm d}\log(f)$ in the special case
$G=(0,\infty)$, $\delta f$ is called the {\it logarithmic derivative}
of $f$.

\begin{Theorem}[Cartan]\label{tslT} Every Maurer--Cartan form $\omega$ on a
simply-connected manifold $\Sigma $ is the logarithmic derivative of
some smooth map $f \colon \Sigma \rightarrow G$. If $f' \colon \Sigma
\rightarrow G$ is a second map with logarithmic derivative $\omega $,
then there exists a unique $g \in G$ such that $f'(x)=f(x)g$.
\end{Theorem}

One says that $f$ is a {\it primitive} of $\omega $. If $\Sigma $ is
only {\itshape connected}, then the obstruction to the existence of a
primitive is known as the {\it monodromy}. Anticipating our later
generalization, we recall two forms of the monodromy here. For further
details see, e.g., \cite[Theorem~7.14, p.~124]{Sharpe_97}.

The {\it global form} of the monodromy is a groupoid morphism
\begin{equation}
 \Omega \colon\ \Pi(\Sigma) \rightarrow G,\label{otter}
\end{equation}
where $\Pi(\Sigma) $ is the fundamental groupoid $\Pi(\Sigma )$ of
$\Sigma $. By definition, an element of $[\gamma] \in \Pi(\Sigma )$ is
the homotopy equivalence class of a path
$\gamma \colon [0,1] \rightarrow \Sigma $ (endpoints fixed). Since the
interval $[0,1]$ is simply-connected, the Maurer--Cartan form
$\gamma^*\omega $ on $[0,1]$ admits, by Cartan's theorem, a unique
primitive $\Gamma \colon [0,1] \rightarrow G$ satisfying
$\Gamma (0)=1_G$, known as the {\it development} of $\omega $ along
the path $\gamma $. One shows that $\Gamma (1) $ depends only on the
class $[\gamma]$ and one defines $\Omega([\gamma])=\Gamma (1)$.

If $\omega $ is the logarithmic derivative of some map $f \colon
\Sigma \rightarrow G$, then
$\Omega([\gamma])=f(\gamma(1))f(\gamma(0))^{-1}$. In~particular,
fixing some $x_0 \in \Sigma $,
\begin{equation}
f(x)=\Omega([\gamma])f(x_0),
\label{bangee}
\end{equation}
where $\gamma \colon [0,1] \rightarrow \Sigma $ is any
path from $x_0$ to $x$. If $\omega $ is an arbitrary Maurer--Cartan
form, then we attempt to {\em define} a primitive $f \colon \Sigma
\rightarrow G$ by \eqref{bangee}. The group of all elements of
$\Pi(\Sigma)$ beginning and ending at $x_0$ is the fundamental group
$\pi_1(\Sigma , x_0)$ and $f$ is well-defined if the restriction of
$\Omega$ to a group homomorphism $\Omega_{x_0} \colon
\pi_1(\Sigma,x_0) \rightarrow G$~--- which we call the {\it pointed
form} of the monodromy~--- is trivial, i.e., takes on the constant
value $1_G$. This condition is evidently independent of the choice of
fixed point $x_0$.

\subsection*{Complete invariants without lifts}
Global lifts as
described above do not exist in general and Cartan's method has been
largely limited to the {\em local} reconstruction of smooth maps into
a homogeneous space, and the special case of curves ($\dim \Sigma =
1$). This is despite the fact that Theorem \ref{tslT} and the monodromy
obstruction are global results!

A generalization of Theorem \ref{tslT} to smooth maps $f \colon \Sigma
\rightarrow G/H$ obviates the need for lifts. Specifically, what we
present here is a characterization of smooth maps $f \colon \Sigma
\rightarrow M$, where $M$ is an arbitrary space on which some Lie
group $G$ is acting transitively, a subtle but significant change in
viewpoint, as we shall explain in Section \ref{infch}. Our results are
naturally formulated in the language of Lie algebroids, and the proof
is an application of Cartan's fundamental theorems for Lie groups,
known as Lie I, Lie II and Lie III, generalized to Lie groupoids, with
which we will assume some familiarity (see, e.g.,
\cite{Crainic_Fernandes_03,Crainic_Fernandes_11}). Standard
introductions to Lie groupoids and algebroids are
\cite{CannasdaSilva_Weinstein_99,Crainic_Fernandes_11,Dufour_Nguyen_05,Mackenzie_05}.

\subsection*{Logarithmic derivatives}

In Lie algebroid language, a Maurer--Cartan form on~$\Sigma $ is
nothing more than a morphism
$\omega \colon {T \Sigma } \rightarrow {\mathfrak g}$ of Lie
algebroids, and Theorem \ref{tslT} a special case of Lie II, as is
well-known. In~the general setting, we replace the Maurer--Cartan form
on $G$ with the action algebroid ${\mathfrak g} \times M$ asso\-ci\-a\-ted
with the action of $G$ on $M$, and use $f \colon \Sigma \rightarrow M$
to pull ${\mathfrak g} \times M$ back to a~Lie algebroid~$A(f)$ over
$\Sigma $. Of course this pullback must be performed {\em in the
 category of Lie algebroids} rather than vector bundles (see, e.g.,
\cite[Section~4.2]{Mackenzie_05}). The composite
$ \delta f \colon A(f) \rightarrow {\mathfrak g}$ of the natural maps
$A(f) \rightarrow {\mathfrak g} \times M \rightarrow {\mathfrak g} $
is a Lie algebroid morphism, which becomes the {\it logarithmic
 derivative} of $f$. The results to be described here show that
$\delta f$~--- or more precisely an appropriate equivalence class of
$\delta f$, see Section~\ref{infch}~--- is a complete invariant of~$f$.

\begin{Example}[logarithmic derivative of an embedding] Suppose
$\Sigma \subset M$ is a submanifold and $f \colon \Sigma \rightarrow
M$ the embedding. Then $A(f)$ is the subbundle of the trivial bundle
${\mathfrak g} \times \Sigma \rightarrow \Sigma$ consisting of all
pairs $(\xi,x)$ having the property that the integral curve on $M$
through $x \in \Sigma $ of the infinitesimal generator $\xi^\dagger$
of $\xi \in {\mathfrak g} $ is tangent to $\Sigma \subset M$ at
$x$. The anchor of $A(f)$ is $(\xi,x) \mapsto \xi^\dagger(x)$ and the
bracket well-defined by
\begin{equation*} [X,Y]=\nabla_{\#X}Y-\nabla_{\#Y}X +\{X,Y\}.
\end{equation*} Here $\nabla $ is the canonical flat connection on
${\mathfrak g} \times \Sigma$ and, viewing sections of ${\mathfrak g}
\times \Sigma $ as ${\mathfrak g} $-valued functions,
$\{X,Y\}(x):=[X(x),Y(x)]_{\mathfrak g} $. The logarithmic derivative
$\delta f$ is the composite $A(f)\hookrightarrow {\mathfrak g} \times
\Sigma \rightarrow {\mathfrak g}$.
\end{Example}

\subsection*{Generalized Maurer--Cartan forms}

Again let $M$ be a smooth manifold on which some Lie group $G$ is
acting from the left transitively~--- what is hereafter referred to as
a {\it homogeneous $G$-space}. With this data fixed, our next task is
to introduce axioms for Lie algebroid morphisms
$\omega \colon A \rightarrow {\mathfrak g}$, where $A$ is a Lie
algebroid over some manifold $\Sigma$, modeled on local properties of
logarithmic derivatives $\delta f$ of smooth maps
$f \colon \Sigma \rightarrow M$.

To this end, observe that {\it logarithmic derivatives map Lie
algebroid isotropy algebras isomorphically onto isotropy algebras of
the action of $G$ on $M$}. Specifically, if we denote the kernel of an
anchor map $\#\colon A \rightarrow T \Sigma$ by $A^\circ $, then, for
any $x \in \Sigma$, $A(f )_x^\circ$ and ${\mathfrak g}_{f(x)}$ have the
same dimension, and
\begin{equation}
\delta f(A(f)_x^\circ) = {\mathfrak g}_{f(x)}.
 \label{uht}
\end{equation}
The following axioms, then, are no stronger than
properties already satisfied by logarithmic derivatives:
\begin{enumerate}\itemsep=0pt
\item[M1] $A$ is transitive.
\item[M2] For some point $x_0 \in \Sigma $, the restriction
$\omega \colon A^\circ_{x_0} \rightarrow {\mathfrak g}$ is injective.
\item[M3] For some (possibly different) point $x_0 \in \Sigma $,
there exists $m_0 \in M$ such that $\omega\big(A^\circ_{x_0}\big) \subset
{\mathfrak g}_{m_0}$, which will be written $x_0\xrightarrow{\omega}m_0$.
\end{enumerate}
Using the shorthand defined in M3, \eqref{uht} reads $x
\xrightarrow{\delta f} f(x)$. A Lie algebroid morphism $\omega \colon
A \rightarrow {\mathfrak g}$ is called a {\it generalized
Maurer--Cartan form} if it satisfies M1--M3.

Contrary to the group case, a Lie algebroid is not necessarily
integrable (the Lie algebroid of some Lie groupoid) and obstructions
to integrability are subtle. See \cite{Crainic_Fernandes_11} for a
fuller discussion and examples. Nevertheless, we have:
\begin{Proposition}\label{frogP}
Assume {\normalfont M1} and {\normalfont M2} hold
and that $\Sigma$ is connected. Then:
 \begin{enumerate}\itemsep=0pt
 \item[$1.$] $A$ is an integrable Lie algebroid.
 \item[$2.$] {\normalfont M2} holds with $x_0$ replaced by an arbitrary
point $x \in \Sigma $.
 \item[$3.$] If {\normalfont M3} holds, then it holds with $x_0$
replaced by an arbitrary point $x \in \Sigma $ $($and suitable choice of
replacement $m \in M$ for $m_0)$.
 \end{enumerate}
\end{Proposition}
We shall see the first assertion readily implies the others. A simple
proof of (1) is not known to us.\footnote{In an earlier
 version of this manuscript the main existence theorem was proven
 without assuming integrability, using substantially more complicated
 arguments, and integrability established post facto.} In Section
\ref{newsec}, where the proposition is proven, we will easily deduce
integrability from Crainic and Fernandes' generalization of Lie III
\cite{Crainic_Fernandes_03}.

\subsection*{Principal primitives}

In fact, the most na\"ive notion of a primitive is not unique ``up to
symmetry''. However, the na\"ive notion will play a role and be given a
name:
\begin{Definition}
 A smooth map $f \colon \Sigma \rightarrow M$ is a {\it principal
 primitive} of a generalized Maurer--Cartan form
 $\omega \colon A \rightarrow {\mathfrak g} $ if there exists a Lie
 algebroid morphism $L \colon A \rightarrow A(f)$ such that the
 following diagram commutes:
 \begin{equation}
 \begin{tikzcd} A \arrow{d}{L} \arrow{r}{\omega } & {\mathfrak g}
 \\ A(f). \arrow[swap]{ru}{\delta f} &
 \end{tikzcd}
 \label{carpet:court}
 \end{equation}
\end{Definition}
 Note that we do not assume $L $ is an isomorphism, or even
that $A$ and $A(f)$ have the same rank.

\subsection*{Monodromy obstructions to the existence of
 primitives}

We now offer this paper's main construction, and formulate the
existence part of our results. Let~$M$ be a homogeneous $G$-space and
$\omega \colon A \rightarrow {\mathfrak g} $ an associated generalized
Maurer--Cartan form. We~are going to explicitly describe the
obstruction to the existence of a principal primitive
$f \colon \Sigma \rightarrow M $ of $\omega $, where $\Sigma $ is the
base of $A$. The most natural description is in terms of some abstract
transitive Lie groupoid ${\mathcal G}$ integrating $A$, whose
existence is guaranteed by Proposition \ref{frogP}, and which we may
take to be source-simply connected, on account of Lie I. In~the next
subsection we will offer a more concrete interpretation using a
generalization of Cartan's development along paths.

According to Lie II, $\omega $ is the derivative of a unique Lie
groupoid morphism
\begin{equation}\label{global}
 \Omega \colon\ {\mathcal G} \rightarrow G,
\end{equation}
which we call the {\it global form} of the monodromy of $\omega $,
being the analogue of \eqref{otter}. Continuing the analogy, we
choose $x_0 \in \Sigma $ and $m_0$ such that M2 and M3 hold, and
attempt to define a~principal primitive
$f \colon \Sigma \rightarrow M$, mapping $x_0$ to $m_0$, by
$f(x)=\Omega(p)\cdot m_0 $, where $p \in {\mathcal G} $ is any arrow
from $x$ to $x_0$. In~this ambition we are successful, so long as $f$
is well-defined, i.e., provided
\begin{gather}
\Omega({\mathcal G}_{x_0}) \subset G_{m_0}.\label{newt}
\end{gather}
Here ${\mathcal G}_{x_0} \subset {\mathcal G}$ denotes
the group of all arrows $p \in {\mathcal G} $ beginning and ending at
$x_0$, and $G_{m_0}$ the isotropy at $m_0 $ of the action of $G$ on
$M$.

Now the algebroid isotropy $A^\circ_{x_0}$ is the Lie algebra of
${\mathcal G}_{x_0} $ and, by our hypothesis M3,
$\omega\big(A^\circ_{x_0}\big) \allowbreak\subset {\mathfrak g}_{m_0}$. Therefore,
\begin{equation} \Omega\big({\mathcal G}_{x_0}^\circ\big) \subset G_{m_0},
\label{warbler}
\end{equation}
where ${\mathcal G}_{x_0}^\circ$ is the connected
component of ${\mathcal G}_{x_0}$. Moreover, as the transitive Lie
groupoid ${\mathcal G} $ has simply-connected source-fibres, there is
a natural exact sequence
\begin{equation*}
1 \rightarrow {\mathcal G}_{x_0}^\circ \rightarrow
{\mathcal G}_{x_0}\xrightarrow{\rho} \pi_1(\Sigma, x_0) \rightarrow 1.
\end{equation*} From this and \eqref{warbler} we obtain a map
$\Omega_{x_0}^{m_0} \colon \pi_1(\Sigma ,x_0) \rightarrow M$
well-defined by
\begin{equation*} \Omega_{x_0}^{m_0}(\rho ( p))=\Omega(p) \cdot m_0.
\end{equation*}
We call this the {\it pointed form} of the
monodromy. By construction, our requirement~\eqref{newt} is equivalent
to $\Omega_{x_0}^{m_0}$ taking a constant value (which is necessarily~$m_0$).

\begin{Theorem}[existence and uniqueness of principal primitives]\label{outlineT}
 A Maurer--Cartan form $\omega{\rm :}$ $A \rightarrow {\mathfrak g} $
 over a connected manifold $\Sigma $ admits a principal primitive
 $f \colon \Sigma \rightarrow M$ if and only if the pointed form of
 the monodromy
 $\Omega_{x_0}^{m_0} \colon \pi_1(\Sigma, x_0) \rightarrow M$ is
 constant for some $($and consequently any$)$ choice of $x_0$ and $m_0$
 with $x_0\xrightarrow{\omega}m_0$. In~that case there is a unique
 principal primitive $f$ of~$\omega $ such that $f(x_0)=m_0$.
\end{Theorem}

\begin{proof}The preceding arguments establish the existence of a
primitive, given constant monodromy. Conversely, given the existence
of a primitive $f$ with $m_0=f(x_0)$, one easily establishes constancy
of the monodromy. For example, an elementary observation stated later
as Proposition \ref{uniquenessP} shows that
 \begin{equation}
 f(x) = \Omega (x_0)\cdot m_0\label{runningshot}
 \end{equation}
 for any arrow $p \in {\mathcal G} $ from $x_0$ to $x$ and so, in particular, $m_0=\Omega (p)\cdot m_0$ for any $p \in
{\mathcal G}_{x_0}$. Since~\eqref{runningshot} applies to any
principle primitive with $m_0=f(x_0)$, the last statement of the
theorem also holds.
\end{proof}

\subsection*{Monodromy as development along $\boldsymbol{A}$-paths}

Let $A$ be a Lie algebroid over a connected manifold $\Sigma $. Then a
piece-wise smooth map $a \colon [0,1] \allowbreak\rightarrow A$, covering an
ordinary path $\gamma \colon [0,1] \rightarrow \Sigma $, is called an
{\it $A$-path} if $\# a(t) = \dot \gamma(t)$, for all $t \in
[0,1]$. Here $\# \colon A \rightarrow T \Sigma $ denotes the anchor of
$A$.

Every $A$-path $a$ can be understood as Lie algebroid morphism $\hat a
\colon {T}I \rightarrow A$ defined by $\hat a (\partial/\partial t)=a$
and all such morphisms arise from $A$-paths. In~particular, given any
Lie groupoid~${\mathcal G} $ integrating $A$, we may apply Lie II,
obtaining a Lie groupoid morphism $I \times I \rightarrow {\mathcal G}
$. The image of $(0, t)$ under this morphism, denoted
\begin{equation*}
\integral_0^t a \in {\mathcal G},
\end{equation*} is an arrow from $\gamma(0)$ to $\gamma(t)$. If $A$ is
a Lie algebra, and ${\mathcal G} $ a Lie group with Lie algebra $A$, then
$a$ is simply a piece-wise smooth path in the Lie algebra, and the
integral above the usual integral to an element in the group. This
familiar case is the one applying in the proposition below:

\begin{Proposition}
Consider a Maurer--Cartan form $\omega \colon A
 \rightarrow {\mathfrak g} $ as in the preceding theorem, and suppose
 $x_0\xrightarrow{\omega}m_0$, for some $x_0 \in \Sigma $ and
 $m_0 \in M$. Let $[\gamma] \in \pi_1(\Sigma,x_0)$ be given and let
 $a \colon [0,1]\rightarrow A$ be {\em any} $A$-path covering
 $\gamma $. Then the monodromy is given by
 \[
 \Omega_{x_0}^{m_0}([\gamma])=\bigg(\integral_0^1 \omega \circ a\bigg)\cdot m_0.
\]
\end{Proposition}
\begin{proof} The proposition follows immediately from the definition
of $\Omega_{x_0}^{m_0}$ and the following ele\-mentary property of
$A$-paths: Every Lie algebroid morphism $\omega \colon A_1 \rightarrow
A_2$ maps $A_1$-paths to $A_2$-paths, and if $\omega $ is the
derivative of a Lie groupoid morphism $\Omega \colon {\mathcal G}^1
\rightarrow {\mathcal G}^2$, then, for any $A_1$-path $a$,
 \begin{equation*}
 \Omega \bigg(\integral_0^t a\bigg) =
\integral_0^t \omega \circ a, \qquad t \in I.\tag*{\qed}
 \end{equation*}\renewcommand{\qed}{}
\end{proof}

\subsection*{Invariants for subgeometry and Bonnet-type theorems} As
far as we know, Cartan's method of moving frames is the only general
technique for obtaining invariants of a submanifold $\Sigma $ of a
Klein geometry $M \cong G/H$, and for proving theorems which
reconstruct the submanifold from its invariants (up to symmetry). The
fundamental theorem of surfaces (Bonnet theorem) is a prototype for
results of this kind. For the special class of parabolic geometries
($G/H$ a flag manifold) an approach based on tractor bundles is
outlined in \cite{Burstall_Calderbank_04} and successfully applied to
conformal geometry (see also \cite{Burstall_Calderbank_10}). These
authors do not describe the monodromy, however, restricting their
attention to the case of simply-connected submanifolds.

While the logarithmic derivative $\delta f$ introduced here delivers a
complete invariant of an embedding $f \colon \Sigma \hookrightarrow
G/H$, it is usually too abstract to be immediately useful. In~\cite{Blaom_G} we take up the problem of ``deconstructing'' this
invariant, and offer illustrations to concrete geometries.

\subsection*{Bracket convention}
Throughout this article, brackets on
Lie algebras and Lie algebroids are defined using {\em
right}-invariant vector fields.

\section{Integrability}
\label{newsec}

In this section we prove Proposition \ref{frogP} by applying Crainic
and Fernandes' generalization of Lie III \cite{Crainic_Fernandes_03}.

Let $A$ be a Lie algebroid over $\Sigma$. In~\cite{Crainic_Fernandes_03, Crainic_Fernandes_11} the kernel of the
anchor of $A \rightarrow T \Sigma$ is denoted by~${\mathfrak g}
$. However, as this conflicts with our use as the Lie algebra of $G$,
we continue to denote the kernel by $A^\circ$. We otherwise follow the
notation and terminology of \cite{Crainic_Fernandes_11}.

In particular, the Weinstein groupoid of $A$ is denoted ${\mathcal
G}(A)$. An element of ${\mathcal G}(A)$ is a certain equivalence class
of $A$-paths. The obstruction to the existence of a bona fide Lie
groupoid integrating $A$ (that is, to the topological groupoid
${\mathcal G}(A)$ being a Lie groupoid) is measured by the {\it monodromy groups} $\tilde {\mathcal N}_{x_0}(A) $, $x_0 \in \Sigma
$. By definition, $\tilde {\mathcal N}_{x_0}(A) $ is the kernel of the
natural homomorphism
\begin{equation}
{\mathcal G}\big(A^\circ_{x_0}\big)\rightarrow {\mathcal G}(A)^\circ_{x_0} \label{dqp}.
\end{equation} On the right-hand side ${}^\circ$ denotes connected
component. At the level of Lie algebroid paths, this homomorphism is
just inclusion. The object on the left is a Lie group, while that on
the right may only be a topological group. Specializing
\cite{Crainic_Fernandes_03} to the transitive case, we have

\begin{Theorem}
 Assuming its base manifold $\Sigma $ is connected, the Lie algebroid
 $A$ is integrable if and only if
 $\tilde {\mathcal N}_{x_0}(A) \subset {\mathcal G}\big(A^\circ_{x_0}\big)$
 is discrete, for some $x_0 \in \Sigma $.
\end{Theorem}

Now assume M1 and M2 hold. Then the restriction $\omega \colon
A^\circ_{x_0} \rightarrow {\mathfrak g} $ is an injection, integrating
to a homomorphism $\Omega \colon {\mathcal G}\big(A^\circ_{x_0}\big)
\rightarrow G$ {\em of Lie groups}, whose kernel $K_0$ is accordingly
discrete. On the other hand, we may include this homomorphism in the
following commutative diagram, whose vertical arrow is~\eqref{dqp}:
\begin{equation*}
 \begin{tikzcd} {\mathcal G}\big(A^\circ_{x_0}\big) \arrow{d}{}
\arrow{r}{\Omega} & G \\ {\mathcal G}(A)^\circ_{x_0}.
\arrow[swap]{ru}{} &
 \end{tikzcd}
\end{equation*} Here the diagonal map is the restriction of the natural
topological groupoid morphism ${\mathcal G}(A) \rightarrow {\mathcal
G}({\mathfrak g}) = G$, i.e., the map sending the equivalence class of
an $A$-path $a$ to the equivalence class of the ${\mathfrak g}$-path
$\omega \circ a$. Commutativity of the diagram implies the kernel of
the vertical map must lie in~$K_0$, but this kernel is, by definition,
$\tilde {\mathcal N}_{x_0}(A)$. This shows $\tilde {\mathcal
N}_{x_0}(A)$ is a discrete subset of~${\mathcal G}\big(A^\circ_{x_0}\big)$ and
the proof of Proposition~\ref{frogP}(1) now follows from the
theorem above.

Now let $x \in \Sigma $ be arbitrary and let ${\mathcal G} $ be a Lie
groupoid integrating $A$, which is necessarily transitive. Then,
assuming $\Sigma $ is connected, there exists an arrow $p \in
{\mathcal G} $ from $x_0$ to $x$. Conjugation $C_p\colon q \mapsto
pqp^{-1}$ maps the isotropy group ${\mathcal G}_{x_0}$ isomorphically
onto ${\mathcal G}_x$. Differentiating, we get a~Lie algebra
isomorphism ${\rm d}C_p \colon A^\circ_{x_0} \rightarrow A^\circ_x$ and a
commutative diagram
\begin{equation*}
 \begin{CD} A^\circ_{x_0} @>{{\rm d}C_p}>> A^\circ_x \\ @V{\omega}VV
@VV{\omega}V \\ {\mathfrak g} @>{\Adjoint_{\Omega(p)}}>> {\mathfrak g}.
 \end{CD}
\end{equation*}
From this observation parts~(2) and~(3) of Proposition \ref{frogP} immediately follow.

\section{The uniqueness of primitives up to symmetry}\label{infch}
This section establishes the uniqueness of primitives, appropriately
defined, up to symmetry. The central result is Theorem
\ref{uniquenessT}. Combining our uniqueness result with Theorem
\ref{outlineT}, we obtain the existence and uniqueness Theorem
\ref{summaryT}.

\subsection*{Symmetries of a homogeneous space}

For the purposes of constructing a theory with the proper invariance,
we have been regarding our fixed data as a homogeneous
$G$-space. While a choice of point $m_0 \in M$ gives us an
identification $M \cong G/G_{m_0}$, formulations depending on a choice
of base point are to be eschewed.\footnote{Geometries in the real
 world do not come with a preferred choice of base-point. Base-points
 are an artifact of Klein's abstraction of geometry, not an intrinsic
 feature.} This decision has a somewhat unexpected consequence,
anticipated by reconsidering the simplest case.

According to Cartan's theorem, a smooth map $f \colon \Sigma
\rightarrow G$ is uniquely determined by its logarithmic derivative,
up to symmetries of $G$. Here a ``symmetry'' is a right group
translation. However, to obtain an invariant version of Cartan's
result we must broaden both the notion of symmetry and what it means
to be a primitive. To see why, consider a smooth map $f \colon \Sigma
\rightarrow M$, where $ M$ is a smooth manifold on which $G$ is acting
transitively and {\em freely}, so that $M \cong G$, up to choice of
base-point. In~order to drop the right-invariant Maurer--Cartan form
$\omega_G$ on $G$ to a~one-form on $ M$, we must suppose here that $G$
is acting on $ M$ from the {\em left}. For then, fixing $ m_0 \in M$
and defining a $\Phi(g) = g \cdot m_0 $, the diffeomorphism $\Phi
\colon G \rightarrow M$ pushes $\omega_G$ forward to a~one-form
$\omega_M $ on $M$ that is independent of the choice of $ m_0$. But
then $\omega_M$ is {\em not} invariant with respect to the action of
$G$~--- rather it is {\em equivariant}, if we regard $G$ as acting on
${\mathfrak g} $ by adjoint action. In~particular, two smooth maps
$f_1, f_2 \colon \Sigma \rightarrow M$ with $f_2(x)=g \cdot f_1(x)$,
$x \in \Sigma $, have, in~general, {\em different} logarithmic
derivatives: $f_2^*\omega_M = \langle\Adjoint_g,f_1^* \omega_M\rangle
$.

To proceed one defines $f \colon \Sigma \rightarrow M$ to be a {\it primitive} of a one-form $\omega $ on~$\Sigma $ if $f^*\omega_M $ and~$\omega $ agree ``up to adjoint action''. The price one pays for this
relaxed definition is that the logarithmic derivative $f^* \omega_M $
only determines $f$ up to a larger class of symmetries of $M$. Under
an identification $M \cong G$, these symmetries consist of the
diffeomorphisms generated by all right {\em and} left translations.

Symmetries in the general case are formalised as follows:
\begin{Definition}\label{symmetriesD}
 Let $M$ be a homogeneous $G$-space. Then a {\it symmetry} of $M$ is
 any diffeomorphism $\phi \colon M \rightarrow M$ for which there
 exists some $l \in G$ such that
 $\phi(g \cdot m) = lgl^{-1} \cdot \phi(m)$ for all
 $g \in G$, $m \in M$.
\end{Definition}

 The symmetries of $M$ form a Lie group henceforth denoted
$\automorphism(M)$. Evidently, $\automorphism(M)$ contains every left
translation $\phi(m):=k \cdot m$, $k \in G$ (take $l=k$). The
following characterization of symmetries is readily verified:

\begin{Proposition}\label{symmetriesP}
 Fix a point $m_0 \in M$ and identify $M$ with the left coset space
 $G/H$, where $H$ denotes the isotropy at $m_0$. Let $l \in G$ be
 arbitrary and suppose $r \in G$ is in the normaliser of $H$, so that
 there exists a map $\phi \colon G/H \rightarrow G/H$ making the
 following diagram commute:
 \begin{equation*}
 \begin{CD} G @>{g~\mapsto~lgr^{-1}}>> G \\ @VVV @VVV \\ G/H
@>{\phi}>> G/H.
 \end{CD}
 \end{equation*}
 Then $\phi $ is a symmetry of $M$ and all
symmetries of $M$ arise in this way. In~other words, the Lie group
$W:=N_G(H)/H$ acts on the left of $G/H$ according to
\[
rH \cdot gH = gr^{-1}H,
\]
$($an action commuting with the left action of $G)$ and
$\automorphism(M)$ is the Lie group generated by both left
translations and those transformations of $M \cong G/H$ defined by the
action of $W$. Here $N_G(H)$ denotes the normaliser of $H $ in $G$.
\end{Proposition}

In contrast to the special case in which $G$ acts freely,
$\automorphism(M)$ is frequently not much larger than the group $G$ of
left translations, in applications of interest to geometers:
\begin{Examples}
\mbox{}
\begin{enumerate}\itemsep=0pt
\item 
Take $M={\mathbb R}^n$, let $H \subset
{\rm GL}(n,{\mathbb R})$ be any linear Lie group whose fixed point
set is the origin, and let $G \cong H \ltimes {\mathbb R}^n $ be the
group of transformations of ${\mathbb R}^n $ generated by translations
and elements of $H$. Then $N_G(H)=H$ and accordingly
$\automorphism(M)=G$.

\item ({\it Affine geometry}) As special cases of item~(1), we may take
$G={\rm GL}(n,{\mathbb R})$ or $G={\rm SL}(n,{\mathbb R})$ and
obtain the affine and equi-affine geometries, with
$\automorphism(M)=G$.

\item 
({\it Euclidean, elliptic and hyperbolic geometry}) Take
$M$ to be one of Riemannian space forms~${\mathbb R}^n$, ${\mathbb S}^n$ or
${\mathbb H}^n $, and let $G$ be the full group of isometries. Then in
every case it is possible to show that each element of $N_G(H)/H$ has
a representative $r \in N_G(H) $ lying in the centre of $G$, and it
follows easily that $\automorphism(M)=G$.

\item ({\it Special elliptic geometry}) Take $M = {\mathbb S}^n$ but let $G$ be the
group of {\em orientation-preserving} isometries,
${\rm SO}(n+1)$. In~this case a little more work reveals that
 \begin{equation*}
 \automorphism(M)=
 \begin{cases}
{\rm SO}(n+1)&\text{if $n$ is odd},\\
{\rm O}(n+1)&\text{if $n $ is even}.
 \end{cases}
 \end{equation*}
 That is, for even-dimensional spheres, we must add
to $G$ the orientation-reversing isometries to obtain the full
symmetry group.

\item 
Suppose $M$ is a homogeneous $G$-space where $G$ is
compact and connected and has trivial centre, and suppose that the
isotropy subgroup $H$ at some point of $M$ is a maximal torus. Then
$N_G(H)/H$ is the Weyl group, well-known to be {\em finite}.

\item ({\it Parabolic geometries}) For a flag manifold $M$, such as a
conformal sphere or projective space, $G$ is a connected semi-simple
Lie group and the isotropy group $H$ is a parabolic subgroup of
$G$. In~this case also $N_G(H)/H$ is known to be finite.
\end{enumerate}
\end{Examples}

\subsection*{Morphisms between Maurer--Cartan forms}

Henceforth we drop the qualification ``generalized'': All Maurer--Cartan
forms and logarithmic derivatives will be understood in the
generalized sense.

With $G$, $M$ and $\Sigma $ fixed as in the Introduction (under
``Generalized Maurer--Cartan forms'') we collect all associated
Maurer--Cartan forms into the objects of a category. In~this category a
morphism $\omega_1 \rightarrow \omega_2$ between objects
$\omega_1 \colon A_1 \rightarrow {\mathfrak g} $ and
$\omega_2 \colon A_2 \rightarrow {\mathfrak g} $ consists of a Lie
algebroid morphism $\lambda \colon A_1 \rightarrow A_2$ covering the
identity on~$\Sigma $ and an element $l \in G$ such that the following
diagram commutes:
 \begin{equation*}
 \begin{CD} A_1 @>{\omega_1}>> {\mathfrak g} \\ @V{\lambda}VV
@VV{\Adjoint_l}V \\ A_2 @>{\omega_2}>> {\mathfrak g}.
 \end{CD}
 \end{equation*} If $\lambda$ is injective, we will say that $\omega_1
\rightarrow \omega_2$ is {\it monic}. The preceding abstractions are
justified by the following observation (strengthened in special cases
in Theorem~\ref{uniquenessT} below):

\begin{Proposition}\label{morphismsP}
 Let $f_1\colon \Sigma \rightarrow M$ be a smooth map into a
 homogeneous $G$-space $M$ and define a second smooth map
 $f_2 \colon \Sigma \rightarrow M$ by $f_2 = \phi\circ f_1$, for
 some $\phi \in \automorphism(M)$. Then $\delta f_1$ and $\delta f_2$
 are isomorphic in the category of Maurer--Cartan forms.
\end{Proposition}

That is, smooth maps $f_1,f_2 \colon \Sigma \rightarrow M$ agreeing up
to a symmetry of $M$ have isomorphic logarithmic derivatives.

\begin{proof} Supposing $f_2 = \phi \circ f_1$, $\phi \in
\automorphism(M)$, define $l \in G$ as in Definition
\ref{symmetriesD}. Then the map $\Adjoint_l \times \phi $, defined by
 \begin{align*}
 (\xi,x) &\mapsto \big(\!\Adjoint_l \xi , \phi(x)\big),
 \\
{\mathfrak g} \times M &\rightarrow {\mathfrak g} \times M
 \end{align*}
 is a Lie algebroid automorphism of the action
algebroid ${\mathfrak g} \times M$ covering $\phi \colon M \rightarrow
M$. In~particular, the composite $A(f_1) \rightarrow {\mathfrak g}
\times M \xrightarrow{\Adjoint_l \times \phi} {\mathfrak g} \times M$
is a Lie algebroid morphism $J$ sitting in a~com\-mutative diagram
\begin{equation*}
 \begin{CD}
 A(f_1) @>{J}>> {\mathfrak g} \times M \\ @VVV @VVV \\
T \Sigma @>>{Tf_2}> TM.
 \end{CD}
\end{equation*}
The vertical arrows indicate anchor maps. Explicitly, we have
\[
J(X) = \big(\!\Adjoint_l \delta f_1 (X), f_2(\lrcorner X)\big),
\]
where $\lrcorner X \in \Sigma$ denotes the base point of $X$.

As $A(f_2)$ is the pullback of ${\mathfrak g} \times M$ under $f_2$,
we obtain, from the universal property of pullbacks, a unique Lie
algebroid morphism $\lambda \colon A(f_1) \rightarrow A(f_2)$ such
that $J $ is the composite
 \[
 A(f_1) \xrightarrow{\lambda} A(f_2) \xrightarrow{\delta f_2}
 {\mathfrak g} \times M.
 \]

This immediately implies commutativity of the diagram
\begin{equation*}
 \begin{CD} A(f_1) @>{\delta f_1}>> {\mathfrak g} \\
@V{\lambda}VV @VV{\Adjoint_l}V \\ A(f_2) @>{\delta f_2}>> {\mathfrak g}.
 \end{CD}
\end{equation*}
One argues that $\lambda$ is an isomorphism by
replacing $\phi $ with $\phi^{-1}$ and reversing the roles of $f_1$
and~$f_2$.
\end{proof}

\subsection*{Primitives}

A smooth map $f \colon \Sigma \rightarrow M$ will be called a {\it primitive} of the Maurer--Cartan form
$\omega \colon A \rightarrow {\mathfrak g} $ if there exists a
morphism $\omega \rightarrow \delta f$. Evidently, every principal
primitive is a primitive.

\subsection*{Maximal Maurer--Cartan forms}

Note that Axioms M1 and M2, together with Proposition \ref{frogP},
imply the following restrictions on the necessarily constant rank of
$A$, whenever $\omega \colon A \rightarrow {\mathfrak g} $ is a
Maurer--Cartan form:
\begin{equation*}
\dimension {\mathfrak g} -\dimension M\le \rank A
\le \dimension {\mathfrak g} -\dimension M + \dimension \Sigma.
\end{equation*}
We say $\omega $ is {\it maximal} if $A$ has maximal
rank, i.e., if
\[
\dim {\mathfrak g}-\rank A = \dimension M - \dimension \Sigma.
\]
In
this case it follows from M3, Proposition \ref{frogP}, and a dimension
count that
\begin{enumerate}\itemsep=0pt
\item[M3$'$.] For any point $x \in \Sigma $ there exists $m \in M$
such that $\omega(A^\circ_x) = {\mathfrak g}_m$.
\end{enumerate} Logarithmic derivatives and ordinary Maurer--Cartan
forms are always maximal.
\begin{Lemma}
 Every morphism $\omega \rightarrow \delta f$ is monic. In~particular, if $\omega $ is maximal, then
 $\omega \rightarrow \delta f$ is an isomorphism.
\end{Lemma}

\begin{proof} A morphism $\omega \rightarrow \delta f$ consists of a
Lie algebroid morphism $\lambda \colon A \rightarrow A(f)$ covering
the identity on~$\Sigma$, and $l \in G$, such that
 \begin{equation}
 \delta f (\lambda(a))=\Adjoint_l \omega(a), \qquad a \in A.\label{cat}
 \end{equation}
 Suppose $\lambda(a)=0$, $a$ an element of $A $ with
base-point $x \in \Sigma $. Since $\lambda$ is a Lie algebroid
morphism covering the identity, we have $\# a =0$, i.e., $a \in
A^\circ_x$. Since $\omega(a)=0$, by \eqref{cat}, Axiom M2 and
Proposition \ref{frogP} imply $a=0$.
\end{proof}
The existence Theorem \ref{outlineT} has the following corollary (of
which we make no further use):
\begin{Corollary}
 Every Maurer--Cartan form
 $\omega \colon A \rightarrow {\mathfrak g} $ with constant monodromy
 has an extension to a maximal Maurer--Cartan form
 $\omega \colon A' \rightarrow {\mathfrak g} $, for some Lie
 algebroid $A' \supset A$.
\end{Corollary}

\begin{proof} By the existence theorem, $\omega $ admits a principal
 primitive $f \colon\! \Sigma \rightarrow M$. That is, there exists a
 morphism $\lambda \colon\! A \rightarrow A(f)$, injective by the
 lemma, whose logarithmic derivative
 $\delta f \colon\! A(f) \rightarrow {\mathfrak g} $ fits into the
 commutative diagram~\eqref{carpet:court}. The logarithmic
 derivative of $f$ is then a maximal Maurer--Cartan form extending
 $\omega $.
\end{proof}

\subsection*{Uniqueness of primitives}

As usual, suppose $G$ acts transitively on $M$, and let
$G_{m_0}^\circ $ denote the connected component of the isotropy
$G_{m_0}$ at some $m_0 \in M$. Then since $G_{m_0}^\circ$ is {\em
 path}-connected, $N_G\big(G_{m_0}\big)\subset N_G\big(G_{m_0}^\circ\big)$.
\begin{Definition}
 We say the isotropy groups of the $G$ action are {\it weakly
 connected} if for some (and hence any) $m_0 \in M$, we have
 $N_G\big(G_{m_0}\big)=N_G\big(G_{m_0}^\circ\big)$.
\end{Definition}

\begin{Example} If $M $ is one of the Riemannian space forms ${\mathbb
R}^n$, ${\mathbb S}^n$ or ${\mathbb H}^n $, and $G$ is the full group of
isometries, then although the isotropy groups of the action of $G$ on
$M$ are not connected, they {\em are} weakly connected.
\end{Example}
A proof of the following central result appears below.

\begin{Theorem}\label{uniquenessT}
 Suppose the action of $G$ on $M$ has weakly connected isotropy
 groups. Let $f_1, f_2 \colon \Sigma \rightarrow M$ be smooth
 maps. Then there exists an isomorphism $\delta f_1 \cong \delta f_2$
 in the category of Maurer--Cartan forms {\em if and only if} there
 exists $\phi \in \automorphism(M)$ such that
 $f_2 = \phi \circ f_1$.
\end{Theorem}

In contrast to the classical setting (Theorem~\ref{tslT})
there may exist more than one choice of $\phi \in \automorphism(M)$
for which $f_2 = \phi \circ f_1$, even if $G$ acts faithfully on
$M$. For example, consider two {\em constant} maps $f_1$, $f_2$.

Combining the theorem with the lemma above, we obtain:
\begin{Corollary}[uniqueness of primitives]\label{uniquenessC}
 If the action of $G$ on $M$ has weakly connected isotropy groups
 then primitives $f \colon \Sigma \rightarrow M$ of a maximal
 Maurer--Cartan form are unique, up to symmetries of $M$.
\end{Corollary}

A non-maximal Maurer--Cartan form may have distinct primitives not
related by a symmetry:
\begin{Example} Let $G$ be the group of isometries of the plane
$M={\mathbb R}^2 $ with Lie algebra ${\mathfrak g} $ identified with
the Killing fields. Let $x, y \colon {\mathbb R}^2 \rightarrow
{\mathbb R} $ denote the standard coordinate functions and let $
\omega \colon T {\mathbb R} \rightarrow {\mathfrak g} $ be the
generalized Maurer--Cartan form\footnote{Actually $\omega $ is an
ordinary Maurer--Cartan form in this case but we are understanding
primitives as maps into~${\mathbb R}^2$, not maps into $G$!} defined
by
\begin{equation*}
\omega \bigg( \frac{\partial }{\partial t}\bigg)=
-y\frac{\partial }{\partial x }+x \frac{\partial }{\partial y}
\qquad\text{(a constant element of $\mathfrak g$).}
 \end{equation*}
 Then for any $r\ge 0$ the map $f (t)=(r\cos t, r\sin
t)$ is a primitive of $\omega $.
\end{Example}

For the proof of the theorem we need one additional observation:

\begin{Proposition}\label{uniquenessP}
 Suppose $f \colon \Sigma \rightarrow M$ is a principal primitive of
 a Maurer--Cartan form $\omega \colon A \rightarrow {\mathfrak g} $.
 Let $\Omega \colon {\mathcal G} \rightarrow G$ be the global form of
 the monodromy of $\omega $, as defined in \eqref{global}. Then, for
 any $x \in \Sigma $, one has $x\xrightarrow{\omega}f(x)$, and for
 any $x_0\in \Sigma $,
 \begin{equation*} f(x)=\Omega(p)\cdot f(x_0),
 \end{equation*} where $p \in {\mathcal G} $ is any arrow from $x_0$
to $x$.
\end{Proposition}
\begin{proof} For some Lie algebroid morphism $\lambda \colon A
\rightarrow A(f)$, we have a commutative diagram
 \begin{equation*}
 \begin{CD} A @>{\omega }>> {\mathfrak g} \\ @V{\lambda}VV @AAA \\
A(f) @>>> {\mathfrak g} \times M.
 \end{CD}
 \end{equation*}
 Since $\lambda$ covers the identity, the claim
$x\xrightarrow{\omega}f(x)$ follows easily from commutativity and the
definition of the bottom map. Let ${\mathcal G}(f)$ denote the
pullback of the action groupoid $G \times M$ by $f$. Since ${\mathcal
G} $ is source-simply-connected, $\lambda$ is the derivative of a Lie
groupoid morphism $\Lambda \colon {\mathcal G} \rightarrow {\mathcal
G}(f)$ and the following diagram commutes (because the composites
being compared have a source-connected domain and identical
derivatives, by the commutativity of the preceding diagram):
 \begin{equation*}
 \begin{CD} {\mathcal G} @>{\Omega }>> G \\ @V{\Lambda}VV @AAA \\
{\mathcal G}(f) @>>> G\times M.
 \end{CD}
 \end{equation*}
 In particular, if we define $F$ to be the composite
Lie groupoid morphism ${\mathcal G}\xrightarrow{\Lambda}{\mathcal
G}(f)\rightarrow G \times M$, then $F$ covers $f \colon \Sigma
\rightarrow M$ and, by the commutativity,
\begin{equation}
F(p)=(\Omega(p),\alpha(p)),\label{rabbit}
\end{equation} where $\alpha $ denotes source projection. But as $F
\colon {\mathcal G} \rightarrow G \times M$ must respect the target
projections, denoted $\beta$, we also have
$f(\beta(p))=\beta(F(p))$. Now \eqref{rabbit} gives
\begin{equation*} f(\beta(p))=\Omega(p) \cdot \alpha(p),
\end{equation*} which proves the proposition.
\end{proof}
\begin{proof}[Proof of theorem (for $\boldsymbol{\Sigma}$ simply-connected)] That
$\delta f_1$ and $\delta f_2$ must be isomorphic when $f_2 = \phi
\circ f_1 $, $\phi \in \automorphism(M)$, is Proposition
\ref{morphismsP}. Suppose $\delta f_1 \cong\delta f_2$ and assume
initially that $\Sigma $ is simply-connected (needed in the proof of
the lemma below). By definition, there exists $l \in G$ and a Lie
algebroid isomorphism $\lambda \colon A(f_2) \rightarrow A(f_1)$ such
that the following diagram commutes:
\begin{gather}
 \begin{CD} A(f_1) @>{\delta f_1}>> {\mathfrak g} \\ @A{\lambda}AA
@VV{\Adjoint_l}V \\ A(f_2) @>{\delta f_2}>> {\mathfrak g}.\label{cod}
 \end{CD}
\end{gather}
Arbitrarily fixing a point $x_0 \in \Sigma$,
\eqref{uht} gives
\begin{equation}
\delta f_1\big(A(f_1)_{x_0}\big) = {\mathfrak g}_{f_1(x_0)},\qquad
\delta f_2\big(A(f_2)_{x_0}\big) = {\mathfrak g}_{f_2(x_0)}.\label{hws}
\end{equation}

For $i=1$ or $2$, let $\Omega_i \colon {\mathcal G}^i \rightarrow G$
denote the global form of the monodromy of $\delta f_i$, as defined at
\eqref{global}. The Lie algebroid of ${\mathcal G}^i$ is $A(f_i)$ and,
by Lie II for Lie groupoids, there is a unique Lie groupoid
isomorphism $\Lambda \colon {\mathcal G}^2\rightarrow {\mathcal G}^1$
whose derivative is $\lambda$. Taking $\omega = \delta f_i$ in the
preceding proposition, we obtain
\begin{equation} f_1(x) = \Omega_1(p_1) \cdot f_1(x_0),\qquad f_2(x) =
\Omega_2(p_2) \cdot f_2(x_0), \label{pto}
\end{equation} whenever $p_i \in {\mathcal G}^i$ is an arrow from
$x_0$ to $x$. By the commutativity of \eqref{cod}, the Lie groupoid
morphisms $\Omega_2 \colon {\mathcal G}^2 \rightarrow G$ and $p_2
\mapsto l \Omega_1({\Lambda}(p_2))l^{-1}$ have the same
derivative, namely $\delta f_2$, so they must coincide, because
${\mathcal G}^2$ is source-connected:
\begin{equation}
\Omega_2(p_2) = l\Omega_1({\Lambda}(p_2))l^{-1},\qquad
p_2 \in {\mathcal G}^2.\label{ksd}
\end{equation}
Since $\lambda$, and hence $\Lambda$, covers the
identity on~$\Sigma $, $p_1 \in {\mathcal G}^1$ is an arrow from $x_0$
to $x$ if and only if $p_2:=\Lambda(p_1)\in {\mathcal G}^2$ is an
arrow from $x_0$ to $x$. This fact and \eqref{ksd} allow us to rewrite
the second equation in \eqref{pto} as $f_2(x)=l \Omega_1(p_1)l^{-1}
\cdot f_2(x_0)$. Or, choosing $r \in G$ such that
\begin{equation} f_2(x_0)=lr^{-1}\cdot f_1(x_0),\label{five}
\end{equation} we have
\begin{equation} f_1(x) = \Omega_1(p_1) \cdot f_1(x_0),\qquad f_2(x) =
l\Omega_1(p_1)r^{-1} \cdot f_1(x_0), \label{pto2}
\end{equation} whenever $p_1 \in {\mathcal G}^1$ is an arrow from
$x_0$ to $x$.
\begin{Lemma} $r \in G$ lies in the normaliser of $G_{f_1(x_0)}$.
\end{Lemma}

Assuming the lemma holds, there exists, by the
characterization of symmetries in Proposition~\ref{symmetriesP}, an
element $\phi \in \automorphism(M)$ well-defined by $\phi(g \cdot
f_1(x_0))=lgr^{-1}\cdot f_1(x_0)$. Then \eqref{pto2} gives us
$f_2(x)=\phi(f_1(x))$, as required.
\end{proof}
\begin{proof}[Proof of lemma] Since we assume the isotropy groups of
the action of $G$ on $M$ are weakly connected, it suffices to show $r
\in N_G\big(G_{f_1(x_0)}^\circ\big)$. We claim
 \begin{gather}
 \Omega_1\big({\mathcal G}^1_{x_0}\big)=G_{f_1(x_0)}^\circ ,\label{hws2}
 \\
\Omega_2\big({\mathcal G}^2_{x_0}\big)=G_{f_2(x_0)}^\circ.\label{hws3}
 \end{gather}
 Since ${\mathcal G}^i$ is transitive and
source-connected ($i \in \{1,2\}$) the restriction of the target
projection of ${\mathcal G}^i$ to the source-fibre over $x_0$ is a
principal ${\mathcal G}^i_{x_0}$-bundle over $\Sigma $. Since we
assume $\Sigma $ is simply-connected, ${\mathcal G}^i_{x_0}$ is
connected, by the long exact homotopy sequence for this principal
bundle. It~follows that \eqref{hws2} and \eqref{hws3} are consequences
of their infinitesimal analogues, which already appear in \eqref{hws}
above.

 Because $\Lambda \colon {\mathcal G}^2 \rightarrow {\mathcal G}^1$
is a Lie groupoid isomorphism covering the identity, we have
 \begin{equation}
 \Lambda \big({\mathcal G}^2_{x_0}\big) ={\mathcal G}^1_{x_0}.\label{kk}
 \end{equation}
 We now compute
 \begin{align*}
 r G_{f_1(x_0)}^\circ r^{-1} &= l^{-1}G_{lr \cdot
f_1(x_0)}^\circ l = l^{-1}G_{f_2(x_0)}^\circ l=l^{-1} \Omega_2
\big({\mathcal G}^2_{x_0}\big)l = \Omega_1
\big(\Lambda\big({\mathcal G}^2_{x_0}\big)\big)
\\
&=\Omega_1 \big({\mathcal
G}^1_{x_0}\big)=G_{f_1(x_0)}^\circ.
 \end{align*}
 The second and subsequent equalities in this
computation follow from equations \eqref{five}, \eqref{hws3},
\eqref{ksd}, \eqref{kk} and \eqref{hws2} respectively.
\end{proof}

\begin{proof}[Proof of theorem {\normalfont(general case)}]\sloppy If $\delta f_1 \cong
\delta f_2$ but $\Sigma $ is not simply-connected, then \mbox{$\delta(f_1
\circ \pi)\cong \delta(f_2 \circ \pi)$}, where $\pi \colon \tilde
\Sigma \rightarrow \Sigma $ denotes the universal covering map, as it
is not difficult to see. By the result just proven in the
simply-connected case, there exists $\phi \in \automorphism(M)$ such
that $f_1 \circ \pi = \phi \circ f_2 \circ \pi $. But as $\pi $ is
surjective, this immediately implies $f_1 = \phi \circ f_2 $.
\end{proof}

\subsection*{Summary of results}

Suppose $f$ is a primitive of a Maurer--Cartan form $\omega $, so that
$\delta f (\lambda(X))=\Adjoint_l \omega(X)$, for some Lie algebroid
morphism $\lambda $ and element $l \in G$. Then it is not hard to show
that $f'(x)=l^{-1} \cdot f(x)$ defines a {\em principal} primitive of
$\omega $. That is, the existence of primitives already implies the
existence of principal primitives. We may therefore summarise the
results cited in the Introduction and our uniqueness result, Corollary
\ref{uniquenessC}, as follows:

\begin{Theorem}[main theorem]\label{summaryT}
 Let $M$ be a homogeneous $G$-space and
 $\omega \colon A \rightarrow {\mathfrak g} $ an associated
 generalized Maurer--Cartan form, where $A$ is a Lie algebroid over
 some manifold $\Sigma $. Then~$A$ is integrable. Furthermore,
 $\omega $ admits a primitive $f \colon \Sigma \rightarrow M$ if and
 only if it has constant monodromy
 $\Omega_{x_0}^{m_0} \colon \pi_1(\Sigma, x_0) \rightarrow M$, for
 some choice of $x_0 \in \Sigma $ and $m_0 \in M$ with
 $x_0\xrightarrow{\omega}m_0$. Assu\-ming~$\omega $ is maximal, and the isotropy groups of the action of $G$ on $M$ are weakly
 connected, the primitive $f$ is unique up to symmetry.
\end{Theorem}

We reiterate that ``symmetry'' is to be understood in the sense Definition \ref{symmetriesD}.

\subsection*{Acknowledgements} The author is indebted to a referee who
contributed the direct proof of integrability in Section~\ref{newsec}. This substantially simplified the proof of the existence
theorem appearing in earlier manuscripts. We thank Yuri Vyatkin, Sean
Curry, Andreas \v{C}ap, and Rui Fernandes for helpful discussions.

\pdfbookmark[1]{References}{ref}
\LastPageEnding

\end{document}